\documentclass[11pt]{article}
\usepackage{amsmath,amsthm,amssymb,xypic}

\theoremstyle{plain}  

\newtheorem{thm}{Theorem}[section]
\newtheorem{prop}[thm]{Proposition}
\newtheorem{lem}[thm]{Lemma}

\theoremstyle{definition}

\newtheorem{defn}[thm]{Definition}

\theoremstyle{remark}

\newtheorem{rem}[thm]{Remark}

\newcommand{\bA}{\mathbb A}
\newcommand{\bC}{\mathbb C}
\newcommand{\bF}{\mathbb F}
\newcommand{\bN}{\mathbb N}
\newcommand{\bP}{\mathbb P}
\newcommand{\bQ}{\mathbb Q}
\newcommand{\bR}{\mathbb R}
\newcommand{\bZ}{\mathbb Z}
\newcommand{\bG}{\mathbb G}

\newcommand{\cO}{\mathcal O}
\newcommand{\cX}{\mathcal X}
\newcommand{\cY}{\mathcal Y}

\newcommand{\lin}{\text{---}}

\newcommand{\map}{\rightarrow}

\DeclareMathOperator{\Hom}{Hom}

\DeclareMathOperator{\Def}{Def}

\DeclareMathOperator{\Spec}{Spec}

\title{Degenerations of del Pezzo surfaces I}
\author{Paul Hacking and Yuri Prokhorov}
\begin{document}
\maketitle

\section{Introduction} \label{intro}

Let $X$ denote a complex surface with quotient singularities.
We say $X$ is a \emph{Manetti surface} if it admits a smoothing to $\bP^2$.
That is, there exists a family $\cX/(0 \in T)$ over the germ of a curve with 
special fibre $X$ and general fibre $\bP^2$.
These surfaces were first studied in \cite{M}.
We prove the following classification of Manetti surfaces.

\begin{thm}\label{manetti}
The Manetti surfaces are as follows:
\begin{enumerate}
\item The weighted projective planes of the form $\bP(a^2,b^2,c^2)$,
where $(a,b,c)$ is a solution of the Markov equation $a^2+b^2+c^2=3abc$.
\item All $\bQ$-Gorenstein deformations of the surfaces in (1).
\end{enumerate}
\end{thm}

Here, we say $X'$ is a \emph{$\bQ$-Gorenstein deformation} of $X$ if there is a family $\cX/(0 \in T)$ with special fibre
$X$ and general fibre $X'$, such that the total space $\cX$ is $\bQ$-Gorenstein.
The $\bQ$-Gorenstein condition is natural from the point of view of the minimal model program, and simplifies the deformation 
theory considerably.
Note in particular that if $X'$ is a $\bQ$-Gorenstein deformation of $X$ then $K_{X'}^2=K_X^2$ .

Each Manetti surface has no locally trivial deformations, and there is one $\bQ$-Gorenstein deformation parameter 
for each singularity. So, all Manetti surfaces are obtained from a weighted projective plane $\bP(a^2,b^2,c^2)$ as 
above by smoothing some subset of the singularities.

Every solution of the Markov equation is obtained from the solution $(1,1,1)$ by a sequence of ``mutations'' of the form 
$$(a,b,c) \mapsto (a,b,3ab-c).$$
The graph which has vertices labelled by the solutions and edges joining solutions related by a single mutation 
is an infinite tree with vertices of degree $3$.

Given a smooth 3-fold, the minimal model program constructs a birational model $Y$ such that either
$K_Y$ is non-negative, or there is a fibration $Y \rightarrow S$ such that $K_Y$ is negative on fibres.
In the second case, $Y/S$ is called a \emph{Mori fibre space}.
If $Y/S$ is a Mori fibre space with fibre dimension $2$, then the general fibre is a smooth del Pezzo surface and 
we say $Y/S$ is a \emph{del Pezzo fibration}. 
The Manetti surfaces are precisely the surfaces with quotient singularities which occur as fibres of del Pezzo fibrations 
with general fibre $\bP^2$.
More generally, we study the surfaces $X$ with quotient singularities which occur as singular fibres of 
del Pezzo fibrations with arbitrary general fibre.
We give a complete classification in the case $\rho(X)=1$ and $K_X^2 \ge 5$.
We also obtain a list of examples with $\rho>1$ as deformations of the surfaces with $\rho=1$.

\begin{thm}\label{generalised}
Let $X$ be a normal surface with quotient singularities such that 
$-K_X$ is ample and $X$ admits a $\bQ$-Gorenstein smoothing. 
Assume $\rho(X)=1$ and $K_X^2 \ge 5$. Then $X$ is one of the following:
\begin{enumerate}
\item A weighted projective plane of the form
\begin{enumerate}
\item $\bP(a^2,b^2,c^2)$, where $a^2+b^2+c^2=3abc$ ($K^2=9$).
\item $\bP(a^2,b^2,2c^2)$, where $a^2+b^2+2c^2=4abc$ ($K^2=8$).
\item $\bP(a^2,2b^2,3c^2)$, where $a^2+2b^2+3c^2=6abc$ ($K^2=6$).
\item $\bP(a^2,b^2,5c^2)$, where $a^2+b^2+5c^2=5abc$ ($K^2=5$).
\end{enumerate}
\item A $\bQ$-Gorenstein deformation of one of the surfaces in (1).
\end{enumerate}
Conversely, for each such surface $X$, the divisor $-K_X$ is ample and $X$ admits a $\bQ$-Gorenstein smoothing 
to a del Pezzo surface.
\end{thm}

\noindent\emph{Acknowledgements:}
We thank Igor Dolgachev for his helpful comments on a preliminary version of this paper.
We also thank Brendan Hassett and S\'{a}ndor Kov\'{a}cs for organising the conference ``Compact moduli spaces and birational geometry''
where the results of this paper were presented, and the American Institute of Mathematics for hosting the conference and generous support.
The work of the second author was partially supported by Grant no. 05-01-00353-a of the Russian Foundation for Basic Research. 

\section{$T$-singularities} \label{secT}

\begin{defn} \cite[Def.~3.7]{KSB} \label{defT}
Let $P \in X$ be a quotient singularity of dimension $2$.
We say $P \in X $ is a \emph{$T$-singularity} if it admits a $\bQ$-Gorenstein smoothing.
That is, there exists a deformation of $P \in X$ over the germ 
of a curve such that the total space is $\bQ$-Gorenstein and the general fibre is smooth.
\end{defn}

For  $n,a \in \bN$ with $(a,n)=1$, let $\frac{1}{n}(1,a)$ denote the cyclic quotient singularity 
$(0 \in \bA^2_{u,v} / \mu_n)$ given by 
$$\mu_n \ni \zeta \colon (u,v) \mapsto (\zeta u, \zeta^a v).$$
Here $\mu_n$ denotes the group of $n$th roots of unity.

\begin{prop}\cite[Prop.~3.10]{KSB}
A $T$-singularity is either a Du~Val singularity or a cyclic quotient singularity of the form 
$\frac{1}{dn^2}(1,dna-1)$ for some $d,n,a \in \bN$ with $(a,n)=1$.
\end{prop}

\begin{defn}
Let $d \in \bN$. 
A \emph{$T_d$-singularity} is a $T$-singularity of the form $\frac{1}{dn^2}(1,dna-1)$.
\end{defn}

We briefly recall the notion of the canonical covering (or index one cover) of a $\bQ$-Gorenstein singularity.
Let $P \in X$ be a normal surface singularity such that $K_X$ is $\bQ$-Cartier. Let $n \in \bN$ be
the least integer such that $nK_X$ is Cartier, the \emph{index} of $P \in X$.
The \emph{canonical covering} $p : Z \rightarrow X$ of $P \in X$ is a Galois cover of $X$ with group $\mu_n$, 
such that $Z$ is Gorenstein and $p$ is \'{e}tale in codimension $1$. Explicitly, 
$$Z= \underline{\Spec}_X(\cO_X \oplus \cO_X(K_X) \oplus \cdots \oplus \cO_X((n-1)K_X))$$ 
where the multiplication in $\cO_Z$ is given by fixing an isomorphism $\cO_X(nK_X) \cong \cO_X$.
The canonical covering is uniquely determined up to isomorphism (assuming that the ground field $k=\bC$ and we work locally 
analytically at $P \in X$).

The singularity $\frac{1}{dn^2}(1,dna-1)$ has index $n$ and canonical covering $\frac{1}{dn}(1,-1)$, 
the Du~Val singularity of type $A_{dn-1}$. We have an identification
$$\frac{1}{dn}(1,-1)=((xy=z^{dn}) \subset \bA^3_{x,y,z}),$$ 
where $x=u^{dn}$, $y=v^{dn}$, and $z=uv$. Passing to the quotient gives
$$\frac{1}{dn^2}(1,dna-1) = \left( (xy=z^{dn}) \subset \frac{1}{n}(1,-1,a) \right).$$ 
A $\bQ$-Gorenstein smoothing of $\frac{1}{dn^2}(1,dna-1)$ is given by 
$$(xy=z^{dn}+t) \subset \frac{1}{n}(1,-1,a) \times \bA^1_t.$$

Let $P \in X$ be a $T$-singularity.
A \emph{$\bQ$-Gorenstein deformation} $(P \in \cX)/(0 \in S)$ of $P \in X$ is by definition a deformation induced by 
an equivariant deformation of the canonical covering \cite[Def.~3.1]{H}.
If $(0 \in S)$ is a smooth curve this is equivalent to requiring that $\cX$ is $\bQ$-Gorenstein.
A singularity of type $\frac{1}{dn^2}(1,dna-1)$ has versal $\bQ$-Gorenstein deformation
$$(xy=z^{dn}+a_{d-1}z^{(d-1)n}+\cdots+a_0) \subset \frac{1}{n}(1,-1,a) \times \bA^d,$$
where $a_0,\ldots,a_{d-1}$ are the coordinates on the base $(0 \in \bA^d)$. In particular, 
$T$-singularities have unobstructed $\bQ$-Gorenstein deformations.
(Note of course that the Du~Val singularities are hypersurface singularities, so in particular Gorenstein.
Thus the $\bQ$-Gorenstein deformations are just the ordinary deformations, which are unobstructed).

We will need the following observation.

\begin{lem} \label{Tdeformation}
Let $(P \in \cX)/(0 \in T)$ be a $\bQ$-Gorenstein deformation of $\frac{1}{dn^2}(1,dna-1)$ over the germ of a curve.
Then the possible singularities of the general fibre are as follows: 
either $A_{e_1-1},\ldots,A_{e_s-1}$ or \mbox{$\frac{1}{e_1n^2}(1,e_1na-1)$}, $A_{e_2-1},\ldots, A_{e_s-1}$,
where $e_1,\ldots,e_s$ is a partition of $d$.
\end{lem}
\begin{proof}
The general fibre has the form 
$$(xy=z^{dn}+a_{d-1}z^{(d-1)n}+\cdots+a_0) \subset \frac{1}{n}(1,-1,a),$$
some $a_0,\ldots,a_{d-1} \in k$.
Write $z^{dn}+a_{d-1}z^{(d-1)n}+\cdots+a_0= \prod (z^n-\gamma_i)^{e_i}$ where the $\gamma_i$ are distinct.
The result follows (the second case occurs when $\gamma_i = 0$ for some $i$).
\end{proof}

\section{Deformations}

We prove (Prop.~\ref{unobs}) that a del Pezzo surface with $T$-singularities has unobstructed $\bQ$-Gorenstein
deformations. We also give (Prop.~\ref{K^2}) a formula relating $K_X^2$ and $\rho(X)$ for such a surface $X$.

\begin{lem} \label{RR}
Let $X$ be a surface with $T$-singularities such that $-K_X$ is ample. 
Then $h^0(\cO_X(-K_X))=1+K_X^2$.
\end{lem}

\begin{proof}
For $X$ a projective surface with quotient singularities and $D$ a divisor on $X$,
we have the singular Riemann--Roch formula \cite{B},\cite{R}
$$\chi(\cO_X(D))=\chi(\cO_X)+ \frac{1}{2}D(D-K_X) + \sum c_P(D),$$
where the $c_P(D)$ are the correction terms corresponding to the points 
$P \in X$ where the divisor $D$ is not Cartier.
We prove below that $c_P(-K_X)=0$ for $P \in X$ a $T$-singularity.
So $\chi(\cO_X(-K_X))=\chi(\cO_X)+K_X^2$.
Since $-K_X$ is ample, we have $H^i(\cO_X(-K_X))=H^i(\cO_X)=0$ for $i>0$ by Kawamata--Viehweg vanishing. 
Hence $h^0(\cO_X(-K_X))=1+K_X^2$, as required.

Let $P \in X$ be a $T$-singularity.
There exists a projective surface $Y$ with a single singularity $Q \in Y$ isomorphic to $P \in X$
and a $\bQ$-Gorenstein smoothing $\cY/(0 \in T)$ of $Y$ over the germ of a curve
(this is a special case of the globalisation theorem \cite[App.]{L}).
By comparing Riemann--Roch on $Y$ and the general fibre of $\cY/T$, we deduce that $c_P(-K_X)=0$, 
as required. Alternatively, the correction term $c_P(-K_X)$ may be computed using a resolution of $X$, 
cf. \cite[p.~312, Pf. of Thm.~1.2]{B}. 
\end{proof}

\begin{prop} \label{unobs}
Let $X$ be a surface with $T$-singularities such that $-K_X$ is ample.
Then there are no local-to-global obstructions for deformations of $X$.
Hence $X$ has unobstructed $\bQ$-Gorenstein deformations, and 
$X$ admits a $\bQ$-Gorenstein smoothing.
\end{prop}

\begin{proof}(cf. \cite{M}, p.~113, Proof of Thm.~21)
The local-to-global obstructions lie in $H^2(T_X)$. 
We have $H^2(T_X)=\Hom(T_X,\cO_X(K_X))^*$ by Serre duality.
Since $H^0(\cO_X(-K_X)) \neq 0$ by Lemma~\ref{RR}, we have an inclusion
$$\Hom(T_X,\cO_X(K_X)) \subset \Hom(T_X,\cO_X) = H^0(\Omega_X^{\vee\vee}).$$
Here $\Omega_X^{\vee\vee}$ is the double dual or reflexive hull of $\Omega_X$.
Let $\pi \colon \tilde{X} \rightarrow X$ be the minimal resolution of $X$.
Then $\pi_*\Omega_{\tilde{X}} = \Omega_X^{\vee\vee}$ by \cite[Lem.~1.11]{S} because $X$ has only 
quotient singularities, and $h^0(\Omega_{\tilde{X}})=0$ because $\tilde{X}$ 
is rational. Combining, we deduce $H^2(T_X)=0$. It follows that $X$ has unobstructed $\bQ$-Gorenstein
deformations since $T$-singularities have unobstructed $\bQ$-Gorenstein deformations (see Sec.~\ref{secT}).
\end{proof}

\begin{defn}
For $P \in X$ a $T$-singularity, 
let $\mu_P$ denote the Milnor number of a $\bQ$-Gorenstein smoothing of $P \in X$.
\end{defn}

\begin{rem}
The number $\mu_P$ is well defined because the $\bQ$-Gorenstein deformation space of a $T$-singularity
is smooth.
\end{rem} 

If $P \in X$ is a Du Val singularity of type $A_n$, $D_n$ or $E_n$ then $\mu_P=n$. 
If $P \in X$ is of type $\frac{1}{dn^2}(1,dna-1)$ then $\mu_P=d-1$ \cite[Sec.~3]{M}.

\begin{prop} \label{K^2}
Let $X$ be a rational surface with $T$-singularities.
Then 
$$K_X^2+\rho(X)+\sum_P \mu_P = 10$$
where the sum is over the singular points $P \in X$.
\end{prop}
\begin{proof}
Let $\tilde{X} \rightarrow X$ be the minimal resolution of $X$.
Then $K_{\tilde{X}}^2+\rho(\tilde{X})=10$ since $\tilde{X}$ is rational.
An explicit computation shows that the change in $K^2+\rho$ 
due to the singularity $P \in X$ equals $\mu_P$ (see, e.g., \cite[p.~111]{M}).
Alternatively, this follows from Noether's formula applied to a smoothing, cf. \cite[4.1f]{L}.
\end{proof}

\section{The Classification}

\begin{thm} \label{toric}
Let $X$ be a toric surface with $T$-singularities such that $\rho(X)=1$ and 
$K_X^2 \ge 5$.
Then $X$ is one of the following surfaces.
\begin{enumerate}
\item $\bP(a^2,b^2,c^2)$, where $a^2+b^2+c^2=3abc$ ($K^2=9$).
\item $\bP(a^2,b^2,2c^2)$, where $a^2+b^2+2c^2=4abc$ ($K^2=8$).
\item $\bP(a^2,2b^2,3c^2)$, where $a^2+2b^2+3c^2=6abc$ ($K^2=6$).
\item $\bP(a^2,b^2,5c^2)$, where $a^2+b^2+5c^2=5abc$ ($K^2=5$).
\end{enumerate}
\end{thm}

The solutions of the equations in Theorem~\ref{toric},(1)--(4) may be described as follows \cite[3.7]{KN}. 
We say a solution $(a,b,c)$ is \emph{minimal} if $a+b+c$ is minimal.
The equations (1),(2),(3) have a unique minimal solution $(1,1,1)$, and (4) has minimal solutions
$(1,2,1)$ and $(2,1,1)$. 
Given one solution, we obtain another by regarding the equation as a quadratic in one of the variables, $c$ (say), 
and replacing $c$ by the other root. 
Explicitly, if the equation is $\alpha a^2 +\beta b^2 + \gamma c^2 =\lambda abc$, then
$$(a,b,c) \mapsto (a,b, \frac{\lambda}{\gamma}ab-c).$$
This process is called a \emph{mutation}.
Every solution is obtained from a minimal solution by a sequence of mutations.

For each equation, we define an infinite graph $\Gamma$ whose vertices are labelled by the solutions, such that two vertices 
are joined by an edge if they are related by a mutation. For equation (1), $\Gamma$ is an infinite tree such that
each vertex has degree $3$, and there is an action of $S_3$ on $\Gamma$ given by permuting the variables $a,b,c$.
The other cases are rather similar, see \cite[3.8]{KN} for details.

\begin{thm} \label{red2max}
Let $X$ be a surface with $T$-singularities such that 
$-K_X$ is ample. Assume that $\rho(X)=1$ and $K_X^2 \ge 5$.
Then $X$ is a $\bQ$-Gorenstein deformation of one of the toric surfaces 
of Thm.~\ref{toric}.

Conversely, if $X$ is an arbitrary $\bQ$-Gorenstein deformation of one of the toric 
surfaces, then $X$ has $T$-singularities and $-K_X$ is ample.
We have $\rho(X)=1$ if and only if the induced deformation of each singularity $P \in X$ 
is either trivial or a deformation of a $T_d$-singularity to an $A_{d-1}$ singularity for some $d \in \bN$.
\end{thm}

\begin{proof}[Proof of Theorem~\ref{toric}]
We first show that $X$ is a weighted projective plane.
A projective toric surface $X$ is given by a complete fan $\Sigma$ in $N_{\bR} \cong \bR^2$, 
where $N \cong \bZ^2$ is the lattice of 1-parameter subgroups of the torus. 
If $\rho(X)=1$ then $\Sigma$ has 3 rays. 
Let $v_0,v_1,v_2 \in N$ be the minimal generators of the rays.
There is a unique relation $a_0v_0+a_1v_1+a_2v_2=0$ where $a_0,a_1,a_2 \in \bN$ are pairwise coprime.
If $v_0,v_1,v_2$ generate the lattice $N$, then $X$ is the weighted projective plane
$\bP(a_0,a_1,a_2)$. In general, let $N'$ denote the lattice generated by $v_0,v_1,v_2$.
Then there is a finite toric morphism $X'=\bP(a_0,a_1,a_2) \rightarrow X$ corresponding to the inclusion 
$N' \subset N$. It has degree $d=|N/N'|$ and is \'{e}tale in codimension $1$.
In particular, we have $K_{X'}^2= d \cdot K_X^2$.
In our situation, $K_X^2 \ge 5$ by assumption. 
Moreover $X'$ has $T$-singularities because a cover of a $T$-singularity which is \'{e}tale in codimension $1$
is again a $T$-singularity (this follows easily from the classification of $T$-singularities).
So $K_{X'}^2 \le 9$ by Proposition~\ref{K^2}. We deduce that $d=1$ and $X=\bP(a_0,a_1,a_2)$.  

The surface $X=\bP(a_0,a_1,a_2)$ has singularities 
$\frac{1}{a_0}(a_1,a_2)$, $\frac{1}{a_1}(a_0,a_2)$, $\frac{1}{a_2}(a_0,a_1)$.
By the classification of $T$-singularities, we have $a_i=d_in_i^2$, where $(a_0,a_1+a_2)=d_0n_0$, etc.
Note that $d_0+d_1+d_2+K_X^2=12$ by Proposition~\ref{K^2}.
Now $K_X \sim -(a_0+a_1+a_2)H$, where $H$ is the ample generator of the class group of $X$, and 
$H^2=1/a_0a_1a_2$. Thus
$$d_0n_0^2+d_1n_1^2+d_2n_2^2= \sqrt{d_0d_1d_2K_X^2}n_0n_1n_2$$
Equations of this form which admit integer solutions were classified in \cite[Sec.~3.5]{KN},
in the context of exceptional collections of vector bundles on del~Pezzo surfaces.
This yields the classification above.
It remains to check that each such surface has $T$-singularities. Let the corresponding equation be
$\alpha a^2 +\beta b^2 + \gamma c^2= \lambda abc$,
and consider the singularity $\frac{1}{\gamma c^2}(\alpha a^2, \beta b^2)$. 
We have $\alpha a^2+\beta b^2 \equiv \lambda abc \equiv 0 \mod \gamma c$. 
Moreover, $c$ is coprime to $a$, $b$ and $\frac{\lambda}{\gamma}$ by the inductive description of the solutions 
of the equation above. So, the singularity is of type $\frac{1}{\gamma c^2}(1,\gamma c q -1)$ for some $q$ with $(q,c)=1$, 
as required.
\end{proof}

The proof of Theorem~\ref{red2max} occupies Sections~\ref{minres} and \ref{pfred2max}.
Theorem~\ref{generalised} follows from Theorem~\ref{red2max}, and Theorem~\ref{manetti} follows from
Theorem~\ref{generalised}. Indeed, a Manetti surface $X$ has $-K_X$ ample and $\rho(X)=1$ \cite[p.~90, Main Thm.]{M},
and a smoothing of a normal surface $X$ to the plane is automatically $\bQ$-Gorenstein \cite[p.~96, Corollary~5]{M}.

\section{The configuration of exceptional curves on the minimal resolution} \label{minres}

The proof of Theorem~\ref{red2max} involves an analysis of the configuration of exceptional curves on the 
minimal resolution $\pi \colon \tilde{X} \rightarrow X$ of the surface $X$. In this section we establish some 
preliminary results. Our starting point is the following theorem:

\begin{thm}\cite[Thm.~11, p.~100]{M} \label{horizcurves}
Let $X$ be a surface with rational singularities such that $\rho(X)=1$ and $h^0(-K_X) \ge 6$.
Then $X$ is rational. 
Let $\pi \colon \tilde{X} \rightarrow X$ be the minimal resolution of $X$.
Assuming $X$ is not isomorphic to $\bP^2$,
there exists a birational morphism $\tilde{X} \map \bF_n$ for some $n$.
Let $d$ be the largest such $n$ and $\mu \colon \tilde{X} \map \bF_d$ a birational 
morphism. Let $B'$ denote the strict transform of the negative section $B \subset \bF_d$.
Let $p \colon \tilde{X} \map \bP^1$ be the composition of $\mu$ with the 
projection $\bF_d \map \bP^1$.
Then $d \ge 2$, and the exceptional locus of $\pi$ is the union of $B'$ and 
the components of degenerate fibres of $p$ of self-intersection at most $-2$, i.e.,
$$\mathrm{Ex}(\pi)= B' \cup \bigcup \{ \Gamma \ | \  p_{\star}\Gamma = 0, \, \Gamma^2 \le -2 \}.$$
Moreover, each degenerate fibre of $p$ contains a unique $(-1)$-curve.
\end{thm}

We describe the possible singularities on $X$ in our case.

\begin{prop} \label{cyclic}
Let $X$ be a surface with $T$-singularities such that $-K_X$ is ample, $\rho(X)=1$, and $K_X^2 \ge 5$.
Then $X$ has cyclic quotient singularities.
\end{prop}
\begin{proof}
Since $K_X^2 \ge 5$,  we have $\sum \mu_P \le 4$ by Proposition~\ref{K^2}. Therefore,
either $X$ has only cyclic quotient singularities, or $X$ has a singularity of
type $D_4$ and possibly some cyclic quotient singularities with $\mu_P=0$. 
In the second case, smoothing the cyclic quotient singularities
we obtain a del Pezzo surface $X'$ with $\rho(X')=1$, $K_{X'}^2=5$, and a single singularity 
of type $D_4$. It is well-known that there are no such surfaces (see, e.g., \cite[Theorem~4.3]{AN}).
\end{proof}
 
Given a cyclic quotient singularity $\frac{1}{n}(1,a)$, let $[b_1,\ldots,b_r]$ be the expansion of $n/a$ as a 
Hirzebruch--Jung continued fraction (cf. \cite[p.~46]{F}).
Then the exceptional locus of the minimal resolution of $\frac{1}{n}(1,a)$ is a string of smooth rational curves of 
self-intersections $-b_1,\ldots,-b_r$. 
The strict transforms of the coordinate lines $(u=0)$ and $(v=0)$ intersect the right and left end components
of the exceptional locus respectively.
\begin{rem} \label{reverse}
Note that $[b_r,\ldots,b_1]$ corresponds to the same singularity as $[b_1,\ldots,b_r]$ with the r\^{o}les of the coordinates 
$u$ and $v$ interchanged. Thus, if $[b_1,\ldots,b_r]=n/a$ then $[b_r,\ldots,b_1]=n/a'$ where $a'$ is the inverse of $a$ 
modulo $n$.
\end{rem}

We recall the description of the minimal resolution of the cyclic quotient singularities of type $T$.
Let a \emph{$T_d$-string} be a string $[b_1,\ldots,b_r]$ which corresponds to a $T_d$-singularity.
\begin{prop}\cite[Thm.~17]{M}, cf. \cite[Prop.~3.11]{KSB} \label{exclocusT_d}
\begin{enumerate}
\item $[4]$ is a $T_1$-string and, for $d \ge 2$, $[3,2,\ldots,2,3]$ (where there are $(d-2)$ $2$'s ) 
is a $T_d$-string.
\item If $[b_1,\ldots,b_r]$ is a $T_d$-string, then so are $[b_1+1,b_2,\ldots,b_r,2]$ and $[2,b_1,\ldots,$ $b_r+1]$.
\item For each $d$, all $T_d$-strings are obtained from the example in (1) by iterating the steps in (2).
\end{enumerate}
\end{prop}

Finally, we describe the possible types of the degenerate fibres of the ruling $p: \tilde{X} \rightarrow \bP^1$ in 
Theorem~\ref{horizcurves}.

\begin{defn} 
Let $a,n \in \bN$ with $a<n$ and $(a,n)=1$. We say the fractions $n/a$ and $n/(n-a)$ are \emph{conjugate}. 
\end{defn}

\begin{lem} \label{conjugate} 
If $[b_1,\ldots,b_r]$ and $[c_1,\ldots,c_s]$ are conjugate, then so are $[b_1+1,b_2,\ldots,b_r]$
and $[2,c_1,\ldots,c_s]$. Conversely, every conjugate pair can be constructed from $[2]$,$[2]$ by a sequence of such 
steps. Also, if $[b_1,\ldots,b_r]$ and $[c_1,\ldots,c_s]$ are conjugate then so are $[b_r,\ldots,b_1]$ and $[c_s,\ldots,c_1]$.
\end{lem}
\begin{proof}
If $[b_1,\ldots,b_r]=n/a$ and $[c_1,\ldots,c_s]=n/(n-a)$ then $[b_1+1,b_2,\ldots,b_r]=(n+a)/a$ and
$[2,c_1,\ldots,c_s]=(n+a)/n$. The last part follows immediately from Remark~\ref{reverse}.
\end{proof}

\begin{lem}(cf. \cite[Pf. of Thm.~18]{M}) \label{fibres}
Let $p \colon \tilde{X} \rightarrow C$ be a birational ruling of a smooth surface.
Suppose $f$ is a degenerate fibre of $p$ such that $f$ contains a unique $(-1)$-curve 
$\Gamma$ and $f \setminus \Gamma$ is a disjoint union of strings of rational curves.
Then the dual graph of $f$ is one of the following types:
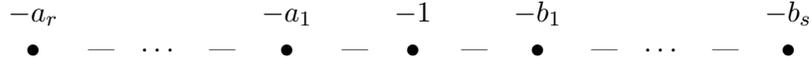
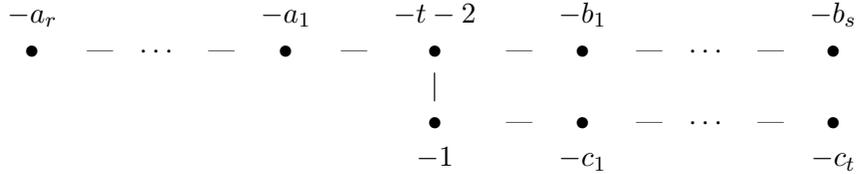
\begin{enumerate}
\item[(I)]
\[
\begin{array}{ccccccccccccc} 
-a_r  &      &        &      & -a_1    &      &  -1     &      &  -b_1   &      &        &      & -b_s    \\
\bullet & \lin & \cdots & \lin & \bullet & \lin & \bullet & \lin & \bullet & \lin & \cdots & \lin & \bullet \\
\end{array}
\]
\item[(II)]
\[
\begin{array}{ccccccccccccc}

-a_r  &      &        &      &  -a_1   &      &  -t-2   &      &   -b_1  &      &        &      & -b_s    \\
\bullet & \lin & \cdots & \lin & \bullet & \lin & \bullet & \lin & \bullet & \lin & \cdots & \lin & \bullet \\
      &      &        &      &         &      &   |     &      &         &      &        &      &         \\
      &      &        &      &         &      & \bullet & \lin & \bullet & \lin & \cdots & \lin & \bullet \\
      &      &        &      &         &      &  -1   &        & -c_1    &      &        &      &  -c_t    \\
\end{array}
\]
\end{enumerate}
In each case $r,s > 0$ and in the second case $t \ge 0$ and $c_1=\cdots=c_t=2$.
Moreover, $[a_1,\cdots,a_r]$ and $[b_1,\cdots,b_s]$ are conjugate. 
Conversely, if these conditions are satisfied then a configuration of smooth rational curves as above is realised as a 
degenerate fibre of a ruled surface.
\end{lem}
\begin{proof}
The ruling $\tilde{X} \rightarrow C$ is obtained from a $\bP^1$-bundle $F \rightarrow C$ by a sequence of blowups.
The lemma follows using the description of conjugate strings in Lemma~\ref{conjugate}.
\end{proof}

\section{Proof of Theorem~\ref{red2max}} \label{pfred2max}

\begin{proof}[Proof of Thm.~\ref{red2max}]
By Lemma~\ref{RR} we have $h^0(-K_X)=1+K_X^2 \ge 6$, so we may apply Theorem~\ref{horizcurves}.
We use the notation of Theorem~\ref{horizcurves}. 
The exceptional locus of the minimal resolution $\pi \colon \tilde{X} \rightarrow X$ is a disjoint union
of strings of rational curves because $X$ has cyclic quotient singularities by Proposition~\ref{cyclic}.
Thus, if $f$ is a degenerate fibre of $p:\tilde{X} \rightarrow \bP^1$, then $f$ is as in Lemma~\ref{fibres}. 
The morphism $\mu \colon \tilde{X} \rightarrow \bF_d$ is an isomorphism over a neighbourhood of the negative section 
$B \subset \bF_d$ because $d$ is maximal by assumption.
So the strict transform $B'$ of $B$ intersects the fibre $f$ in the strict transform $A'$ of the corresponding fibre $A$
of $\bF_d \rightarrow \bP^1$.
The component $A'$ of $f$ is one of the end components labelled $-a_r$ or $-b_s$ in the diagram of Lemma~\ref{fibres}.
We assume for definiteness that it is the left end component labelled $-a_r$.
The connected component of the exceptional locus of $\pi$ containing $B'$ consists of $B'$ together with 
strings of curves meeting $B'$, one for each degenerate fibre $f$. 
Hence there are at most $2$ degenerate fibres, because this connected component is a string of curves.

We show that the surface $X$ is a deformation of a toric surface $Y$ with $T$-singularities such that 
$\rho(Y)=1$. 

We first construct the surface $Y$.
There is a uniquely determined toric blowup $\mu_Y \colon \tilde{Y} \rightarrow \bF_d$ such that
$\mu_Y$ is an isomorphism over the negative section $B \subset \bF_d$, and the degenerate fibres of the ruling 
$p_Y \colon \tilde{Y} \rightarrow \bP^1$ are fibres of type (I) associated to the fibres of $p$ as follows.
Let $f$ be a degenerate fibre of $p$ as in Lemma~\ref{fibres}, and assume that $B'$ intersects the left end component.
If $f$ is of type (I) then the associated fibre $f_Y$ of $p_Y$ has the same form.
If $f$ is of type (II) then $f_Y$ is a fibre of type (I) with self-intersection numbers 
$$-a_r,\ldots,-a_1,-t-2,-b_1,\ldots,-b_s,-1,-d_1,\ldots,-d_u$$ 
Note that the sequence $d_1,\ldots,d_u$ is uniquely determined (see Lemma~\ref{fibres}).
In each case the strict transform $B'$ of $B$ again intersects the left end component of $f_Y$.

Let $Y$ be the toric surface obtained from $\tilde{Y}$ by contracting the strict transform of the negative section of $\bF_d$ 
and the components of the degenerate fibres of the ruling with self-intersection at most $-2$.
For each fibre $f$ of $p$ of type (II) as above, the string of rational curves with self-intersections $-d_1,\ldots,-d_u$
in the associated fibre $f_Y$ of $p_Y$  contracts to a $T_{t+1}$ singularity by Lemma~\ref{T} below.
This singularity replaces the $A_t$ singularity on $X$ obtained by contracting the string of $t$ (-2)-curves in $f$. 
In particular, the surface $Y$ has $T$-singularities. Moreover $\rho(Y)=1$, and $K_Y^2=K_X^2 \ge 5$ by Proposition~\ref{K^2}. 
Thus $Y$ is one of the surfaces listed in Theorem~\ref{toric}. 
A $T_d$-singularity admits a $\bQ$-Gorenstein deformation to an $A_{d-1}$ singularity (see Lemma~\ref{Tdeformation}). 
Hence the singularities of $X$ are a $\bQ$-Gorenstein deformation of the singularities of $Y$.
There are no local-to-global obstructions for deformations of $Y$ by Proposition~\ref{unobs}. 
Hence there is a $\bQ$-Gorenstein deformation $Y'$ of $Y$ with the same singularities as $X$.
We prove below that $X \cong Y'$.
 
Let $f$ be a degenerate fibre of $p$ of type (II) as above and $f_Y$ the associated fibre of $p_Y$.
Let $P \in Y$ be the $T$-singularity obtained by contracting the string 
of rational curves in $f_Y$ with self-intersections $-d_1,\cdots,-d_u$.
Let $\cY/(0 \in T)$ be a $\bQ$-Gorenstein deformation of $Y$ over the germ of a curve
which deforms $P \in Y$ to an $A_{t}$ singularity and is locally trivial elsewhere. Let $Y'$ denote a general fibre.
Let $\hat{\cY} \rightarrow \cY/T$ be the simultaneous minimal resolution of the remaining singularities
(where the deformation is locally trivial), and let $\hat{Y} \rightarrow Y$ be the special fibre and $\hat{Y}' \rightarrow Y'$ 
the general fibre. 
Thus $\hat{Y}$ has a single $T$-singularity and $\hat{Y}'$ a single $A_{t}$ singularity.
The ruling $p_Y \colon \tilde{Y} \rightarrow \bP^1$ descends to a ruling $\hat{Y} \rightarrow \bP^1$; let $A$ be a general fibre of this 
ruling. Then $A$ deforms to a $0$-curve $A'$ in $\hat{Y}'$, which defines a ruling $\hat{Y}' \rightarrow \bP^1$. 
Let $\tilde{Y}' \rightarrow \hat{Y}'$ be the minimal resolution of $\hat{Y}'$ and consider the induced ruling 
$p_{Y'} \colon \tilde{Y}' \rightarrow \bP^1$. Note that the exceptional locus of $\hat{Y} \rightarrow  Y$ deforms without change by 
construction. Moreover, the $(-1)$-curve in the remaining degenerate fibre (if any) of $p_Y$ also deforms. 
There is a unique horizontal curve in the exceptional locus of $\pi_{Y'} : \tilde{Y}' \rightarrow Y'$, and $\rho(Y')=1$ by 
Proposition~\ref{K^2}.
Hence each degenerate fibre of $p_{Y'}$ contains a unique $(-1)$-curve, and the remaining components of the fibre 
are in the exceptional locus of $\pi_{Y'}$. We can now describe the degenerate fibres of $p_{Y'}$. 
If $p_Y$ has a degenerate fibre besides $f_Y$, then $p_{Y'}$ has a degenerate fibre of the same form.
We claim that there is exactly one additional degenerate fibre of $p_{Y'}$, which is of type (II) and has the same form as the fibre 
$f$ of $p$.
Indeed, the union of the remaining degenerate fibres consists of the string of rational curves with self-intersections 
$-a_r,\ldots,-a_1,-t-2,b_1,\ldots,b_s$ (the deformation of the string of the same form in $f_Y$), 
the string of $(-2)$-curves which contracts to the $A_t$ singularity, and some $(-1)$-curves.
The claim follows by the description of degenerate fibres in Lemma~\ref{fibres}.
If there is a second degenerate fibre of $p$ of type (II) we repeat this process.
We obtain a $\bQ$-Gorenstein deformation $Y'$ of $Y$ with minimal resolution $\pi_{Y'} \colon \tilde{Y}' \rightarrow Y'$, 
and a ruling 
$p_{Y'} \colon \tilde{Y}' \rightarrow \bP^1$ such that the exceptional locus of $\pi_{Y'}$ 
has the same form with respect to the ruling $p_{Y'}$ as that of $\pi$ with respect to $p$.

We claim that $X \cong Y'$. Indeed, there is a toric variety $Z$ and, for each fibre $f_i$ of $p$ of type (II), a toric boundary divisor 
$\Delta_i \subset Z$ and points $P_i,Q_i$ in the torus orbit $O_i \subset \Delta_i$, such that $X$ (respectively $Y$) 
is obtained from $Z$ by successively blowing up the points $P_i$ (respectively $Q_i$)  $t_i+1$ times, where $t_i$ is the length of 
the string of $(-2)$-curves in $f_i$. It remains to prove that we may assume $P_i=Q_i$ for each $i$.
Let $T$ be the torus acting on $Z$ and $N$ its lattice of 1-parameter subgroups. Let $\Sigma \subset N_{\bR}$ be the fan corresponding to 
$X$ and $v_i \in N$ the minimal generator of the ray in $\Sigma$ corresponding to $\Delta_i$. 
Then $T_i = (N/\langle v_i \rangle) \otimes \bG_m$ is the quotient torus of $T$ which acts faithfully on $\Delta_i$.
Thus, there is an element $t \in T$ taking $P_i$ to $Q_i$ for each $i$ except in the following case: 
there are two fibres of $p$ of type (II),
and $v_1+v_2=0$. In this case, there is a toric ruling $q \colon Z \rightarrow \bP^1$ given by the projection 
$N \rightarrow N/\langle v_1 \rangle$. The toric boundary of $Z$ decomposes into two sections (given by $\Delta_1,\Delta_2$) and two fibres 
of $q$. But one of these fibres (the one containing $B'$) is a string of rational curves of self-intersections at most $-2$, a contradiction.

The characterisation of the $\bQ$-Gorenstein deformations of the toric surfaces which have Picard number one
is given by Lemma~\ref{Tdeformation} and Proposition~\ref{K^2}.
\end{proof}

\begin{lem}\label{T}
Let $[a_1,\ldots,a_r]$ and $[b_1,\ldots,b_s]$ be conjugate strings.
Then the conjugate of $[a_r,\ldots,a_1,t+2,b_1,\ldots,b_s]$ is a $T_{t+1}$-string.
\end{lem}

\begin{proof}
Let an \emph{$S_t$-string} be a string $[a_r,\ldots,a_1,t+2,b_1,\ldots,b_s]$ as above. Then, by Lemma~\ref{conjugate},
we have
\begin{enumerate}
\item $[2,t+2,2]$ is an $S_t$-string.
\item If $[e_1,\ldots,e_v]$ is an $S_t$-string, then so are $[e_1+1,\cdots,e_v,2]$ and $[2,e_1,\ldots,$ $e_v+1]$.
\item Every $S_t$-string is obtained from the example in (1) by iterating the steps in (2).
\end{enumerate}
The result follows from the description of $T_d$-strings in Proposition~\ref{exclocusT_d} and Lemma~\ref{conjugate}.
\end{proof}

\medskip
\noindent
Paul Hacking, Department of Mathematics, Yale University, PO Box 208283, New Haven, CT 06520, USA; paul.hacking@yale.edu\\
\\
\noindent
Yuri Prokhorov, Department of Higher Algebra, Faculty of Mathematics and Mechanics,  Moscow State Lomonosov University, Vorobievy Gory, 
Moscow, 119 899, RUSSIA; prokhoro@mech.math.msu.su\\

\end{document}